\documentclass{article}
\usepackage{lipsum}  
\usepackage{amsmath}
\usepackage{amssymb}
\usepackage{amsthm}
\usepackage{graphicx}
\usepackage{float}

\providecommand{\keywords}[1]{\textbf{\textit{keywords}} #1}
\providecommand{\MSC}[2]{\textbf{\textit{2000 MSC
}}:26A18,28Dxx,34Cxx,93xx}
\usepackage{blindtext}
\Large{}
\usepackage{subcaption}
\usepackage{hyperref}
\usepackage{geometry} 
\usepackage{lmodern}
\usepackage{algorithm2e}
\graphicspath{{media/}}     
\usepackage{pgfplots}
\pgfplotsset{compat=1.17}
\usepackage{algorithm2e}
\usepackage{lipsum}
\usepackage{mathtools} 
\usepackage{microtype} 
\usepackage{bm} 
\usepackage{hyperref}
\usepackage{url}
\begin{document}
\title{Chaotic Dynamics Derived from the Montgomery Conjecture: Application to Electrical Systems}

\author{
  Zeraoulia Rafik*\\
  Faculty of Material Sciences and Computer Science, Mathematics Department\\
  Khemis Miliana University, Theniet el Had Street, Khemis Miliana (44225), Algeria\\
  Acoustics and Civil Engineering Laboratory\\
  \texttt{zeraoulia@univ-dbkm.dz}\\
  Tel: +213668085320\\
  \and
  Alvaro Humberto Salas\\
  Universidad Nacional de Colombia, Departamento de Matemáticas, Bogotá, Colombia\\
  Research group FIZMAKO, Universidad de Caldas (Colombia)\\
  \texttt{ahsalass@unal.edu.co}\\
  \and
  Ayadi Souad \\
  Faculty of Material Sciences and Computer Science, Physics Department\\
  Khemis Miliana University, Theniet el Had Street, Khemis Miliana (44225), Algeria\\
  Acoustics and Civil Engineering Laboratory\\
  \texttt{souad.ayadi@univ-dbkm.dz}\\
  Tel: +213663447124 / +213775186197
}

\maketitle

\tableofcontents

\

\begin{abstract}
Here, we introduce a novel method for obtaining chaotic dynamics based on the Montgomery conjecture for the pair correlation of zeros of the Riemann zeta function. Motivated by the conjecture, we present a recursive relation that reveals chaotic behavior. Notably, we provide insights into the possible uses of this derived chaotic dynamics in electrical engineering by interpreting it as a unique representation of an electrical system. Furthermore, we investigate the relevance of entropy, bifurcation analysis, and chaos theory in this framework for electrical systems. We look into its applicability to signal processing, stability analysis through bifurcation, and how entropy measures the predictability or unpredictability of electrical signals. Additionally, we discuss the system's strange attractor and its transition to voltage collapse, highlighting the interplay between chaotic dynamics and stability in electrical systems. Furthermore, we analyze the system's energy distribution, taking into account how chaotic dynamics may affect energy allocation or dissipation. Furthermore, we compare the chaotification and Hermiticity of the resulting operators between Yitang dynamics and Montgomery dynamics. To have a better grasp of the spectrum features of each operator, we calculate the eigenvalues for each one obtained from the corresponding dynamics. Our results provide fresh insights into number-theoretic chaotic dynamics and how they might be applied in real-world electrical engineering applications. This work provides encouraging opportunities for further research and technology developments by laying the foundation for creative investigations in system dynamics.

\end{abstract}

\section{Notations Used}\label{sec_notations}

In this section, we provide a list of notations and symbols used throughout the paper for the reader's convenience.

\begin{itemize}
    \item \( F(\omega) \): The function representing the chaotic dynamics in the system.
    \item \( \delta(\omega) \): The Dirac delta function, representing an impulse at \( \omega = 0 \).
    \item \( C \): A constant used in the formulation of the chaotic system.
    \item \( \text{sgn}(\omega) \): The sign function, defined as:
    \[
    \text{sgn}(\omega) = 
    \begin{cases} 
    1 & \text{if } \omega > 0 \\ 
    0 & \text{if } \omega = 0 \\ 
    -1 & \text{if } \omega < 0 
    \end{cases}
    \]
    \item \( \text{sgn}'(\omega) \): The first derivative of the sign function with respect to \( \omega \).
    \item \( \text{sgn}''(\omega) \): The second derivative of the sign function with respect to \( \omega \).
    \item \( r \): A parameter related to the chaotic dynamics.
    \item \( _0F_1\left(; b; z\right) \): The confluent hypergeometric function.
    \item \( _0\tilde{F}_1\left(; b; z\right) \): A variant of the confluent hypergeometric function used in the context of chaotic systems.
    \item \( i \): The imaginary unit, defined as \( i = \sqrt{-1} \).
    \item \( \omega \): The frequency variable in the context of the chaotic dynamics.
\end{itemize}

\section*{Main Results}

\subsection*{1. Derived Electrical System}
The chaotic dynamics inspired by Montgomery's Pair Correlation Conjecture are represented by the following recursive relation, interpreted as an electrical system:
\[
x_{n+1} = 1 - \left( \frac{\sin(\pi/x_n)}{r\pi/x_n} \right)^2 + \frac{R}{L} + \frac{1}{C x_n}
\]
where \(x_n\) represents the system state, and \(R\), \(C\), \(L\), and \(r\) denote the electrical parameters (resistance, capacitance, inductance, and an additional scaling parameter, respectively).

\subsection*{2. Parameter Analysis Insights}
An investigation of how different electrical parameters influence system behavior revealed the following:
\begin{itemize}
    \item \textbf{Inductance Dominance} (\(L > R\)): Leads to smoother transitions and stability in the system.
    \item \textbf{Capacitance and Inductance}: Large values of \(C\) and \(L\) tend to produce periodic behavior, with regular oscillations in the system.
    \item \textbf{Sensitivity to Parameter \(r\)}: Small changes in \(r\) cause significant shifts in system behavior, including the appearance of chaotic dynamics.
\end{itemize}

\subsection*{3. System Dynamics and Analysis}
The derived electrical system demonstrates a range of dynamic behaviors depending on parameter values. Key findings include:
\begin{itemize}
    \item \textbf{Phase Portraits}: Numerical simulations show stable fixed points when \(r = 0.7\), where the system tends to settle into predictable patterns.
    \item \textbf{Limit Cycles and Frequency Analysis}: Fourier analysis reveals the system’s dominant frequencies, corresponding to periodic behavior under certain parameter conditions.
    \item \textbf{Poincaré Section}: Recurrent states appear in the negative region, highlighting the chaotic nature of the system.
\end{itemize}

\subsection*{4. Quantum System Interpretation}
By interpreting the electrical system as a quantum system, the Schrödinger equation with a potential function \(V(x) = A + \frac{1}{C x}\) yields sinusoidal eigenfunctions. These eigenfunctions correspond to discrete energy levels and bound states, indicating parallels between quantum mechanics and signal processing. The eigenfunctions represent fundamental frequency components of the electrical system, essential for stability and performance analysis.

\subsection*{Detailed Energy Dynamics Analysis}

\subsection*{5. Influence of Parameter \(r\)}
The parameter \(r\) plays a crucial role in determining energy storage and distribution in the system. Specifically:
\begin{itemize}
    \item \textbf{Energy Storage Sensitivity}: As \(r\) varies, significant shifts in energy storage behavior occur, affecting both system stability and energy distribution.
    \item \textbf{Energy Variation with Capacitance (\(C\))}: The system shows heightened energy storage in the inductor as \(C\) decreases, suggesting a direct relationship between capacitance and energy dynamics.
\end{itemize}

\subsection*{6. Energy Distribution Implications}
In the context of the derived system represented by the equation:
\[
f(x, r, R, L, C) = 1 - \left( \frac{\sin\left(\frac{\pi}{x}\right)}{r \pi / x} \right)^2 + \frac{R}{L} + \frac{1}{C x}
\]
with parameters \(R = 0.000025\), \(L = 0.00045\), and \(C = 0.73\), the variation of \(r\) has a significant effect on energy storage. Understanding the dynamics of \(r\) could reveal optimal ranges for energy efficiency and system stability.

\subsection*{Eigenvalue and Operator Analysis}

\subsection*{7. Comparison of Eigenvalues}
\begin{itemize}
    \item \textbf{Yitang Dynamics}: The first eigenvalue of the Yitang system, computed using the Implicitly Restarted Arnoldi Method (IRAM), is \(7.11343 - 0.178929i\), indicating complex eigenvalues and potential chaotic behavior.
    \item \textbf{Montgomery Dynamics}: The first eigenvalue of the Montgomery system is \(-2.99007\), a real value, suggesting more regular, predictable dynamics.
\end{itemize}

\subsection*{8. Diagonalizability}
\begin{itemize}
    \item \textbf{Yitang Dynamics}: The complex eigenvalues indicate the possibility of complex eigenvectors, and further analysis is required to determine whether the system is diagonalizable.
    \item \textbf{Montgomery Dynamics}: The real eigenvalues suggest that the Montgomery system's corresponding operator is diagonalizable, with real and linearly independent eigenvectors.
\end{itemize}

\section{Introduction}

Number theory, the study of the properties and relationships of integers, has long been a fertile ground for deep conjectures that continue to intrigue mathematicians. These conjectures often emerge from patterns observed in numerical data, leading to hypotheses whose proof could significantly advance our understanding of fundamental mathematical principles.

A major area of interest in number theory is the distribution of prime numbers. Famous conjectures in this field include the Twin Prime Conjecture, Goldbach’s Conjecture, and the Riemann Hypothesis. Despite their deceptively simple formulations, these conjectures have resisted proof for centuries, captivating mathematicians and driving numerous research efforts.

In recent years, there has been a growing interest in the potential links between number theory and chaotic dynamics. Chaos theory, which studies complex and unpredictable behavior in deterministic systems, has found diverse applications across fields like physics, biology, and economics. The possibility that chaotic behavior could emerge from deterministic number-theoretic systems has inspired investigations into these unexpected connections \cite{3}.

One promising line of research explores how chaotic dynamics can be applied to conjectures in number theory. By analyzing chaotic systems inspired by conjectures such as the Riemann Hypothesis and the Twin Prime Conjecture, researchers seek to uncover new structural insights that could eventually lead to proofs or disproofs.

In this work, we build on our earlier research presented in \cite{24}, where we introduced a novel derivation of chaotic dynamics based on the Montgomery Pair Correlation Conjecture. This conjecture, which describes statistical properties of the non-trivial zeros of the Riemann zeta function \cite{2}, provides a fertile framework for studying the interplay between number theory and chaotic behavior.

Building upon our previous findings, we offer an enhanced analysis and a more comprehensive exploration of chaotic dynamics derived from the Montgomery conjecture. Our investigation not only deepens the understanding of these dynamics but also highlights their potential applications in electrical engineering.

Several recent works have explored advanced control strategies and mathematical techniques that resonate with our study. For instance, Deep et al. \cite{26} introduced a robust moment-matching load frequency control strategy for cyber–physical power systems, addressing communication delays, which has implications for stability analysis. Similarly, Ansari and Raja \cite{27} developed an enhanced cascaded frequency controller optimized by a flow direction algorithm, demonstrating the importance of control optimization in complex dynamical systems. Aryan et al. \cite{28} proposed an equilibrium optimizer tuned frequency-shifted internal model control design for industrial plants, emphasizing structured control in practical applications.

Mukherjee et al. \cite{29} investigated fractional backstepping and Lyapunov strategies for stabilization, which align with the stability concerns inherent in chaotic systems. Kumari et al. \cite{30} explored fuzzy control methods for frequency regulation in power systems, demonstrating the relevance of advanced control techniques in managing nonlinear dynamics. Further, Hu et al. \cite{31} introduced a double predictive proportional-integral control strategy, relevant for handling delays and nonlinearities similar to those found in chaotic number-theoretic models.

Mehta et al. \cite{32} developed a tri-parametric fractional-order controller for integrating systems with time delays, a topic closely linked to our interest in delay-induced chaotic behavior. Mukherjee et al. \cite{33} proposed an optimal fractional Lyapunov-based adaptive control scheme, which has implications for the stability of dynamical systems derived from number-theoretic conjectures. Aryan and Raja \cite{34} extended this line of research with an LFC scheme incorporating renewables and electric vehicle penetration, further highlighting the intersection between applied engineering and mathematical modeling.

Other notable contributions include Aryan et al. \cite{35}, who examined a RIMC-based dual-loop strategy for unstable processes, and their work on internal model control for boiler steam drum regulation \cite{36}, both of which provide insights into structured control methodologies applicable to our chaotic framework. Additionally, Anand et al. \cite{37} investigated type-2 fuzzy controllers optimized via equilibrium algorithms, reinforcing the role of computational optimization in nonlinear control systems.

By drawing parallels between number theory and electrical systems, we uncover new insights that contribute to both fields. This interdisciplinary approach seeks to bridge the gap between chaos theory, number theory, and applied engineering sciences \cite{6}. We aim to foster new avenues for research and innovation at the intersection of these domains \cite{8}.

\section{Montgomery's Pair Correlation Conjecture in Number Theory}

Montgomery's Pair Correlation Conjecture is a significant hypothesis in number theory that explores the statistical distribution of spacings between the non-trivial zeros of the Riemann zeta function. Proposed by Hugh L. Montgomery, this conjecture has established deep connections between the distribution of these spacings and random matrix theory \cite{4}.

\subsection{Background}

The Riemann zeta function, $\zeta(s)$, is a complex-valued function of a complex variable $s$, whose non-trivial zeros are central to understanding the distribution of prime numbers. These non-trivial zeros, also referred to as the Riemann zeros, lie in the critical strip where the real part of $s$ is between 0 and 1, and play a fundamental role in number theory. Montgomery's conjecture focuses on the statistical properties of the gaps between these zeros.

\subsection{The Pair Correlation Conjecture}

Montgomery's Pair Correlation Conjecture predicts that the pair correlation function $R_2(\Delta)$, which measures the probability density of finding two non-trivial zeros with a normalized spacing $\Delta$, converges to a distribution analogous to the eigenvalue spacings of matrices in the Gaussian Unitary Ensemble (GUE) as the imaginary part of the zeros grows large. Specifically, the conjecture can be formulated as:

\begin{equation}
    R_2(\Delta) \to \text{limiting distribution of GUE}, \quad \text{as} \quad \text{Im}(\rho) \to \infty,
\end{equation}

where $\rho$ denotes the non-trivial zeros of the Riemann zeta function, and GUE refers to a class of random matrices whose eigenvalue statistics are well-understood through random matrix theory.

\subsection{The Pair Correlation Function}

The pair correlation function $R_2(\Delta)$ is mathematically defined as:

\begin{equation}
    R_2(\Delta) = \lim_{{T\to\infty \atop T \text{ fixed}}} \frac{1}{T} \sum_{{0 < \text{Im}(\rho_j), \text{Im}(\rho_k) < T \atop |\text{Im}(\rho_j) - \text{Im}(\rho_k)| < \Delta}} 1,
\end{equation}

where $\rho_j$ and $\rho_k$ are the non-trivial zeros of the Riemann zeta function, and $\Delta$ represents the normalized spacing between their imaginary parts.

\subsection{Pair Correlation Formula with Sine Function}

Montgomery's conjecture can also be represented by a pair correlation formula involving the sine function, which approximates the behavior of the zeros. The formula is given by:

\begin{equation}
    R_2(\Delta) = \frac{2}{\pi} \left(\Delta \sin\left(\frac{\pi \Delta}{2}\right) + \frac{\sin^2\left(\frac{\pi \Delta}{2}\right)}{2}\right),
\end{equation}

where $\Delta$ denotes the normalized spacing between the non-trivial zeros. This formula highlights the oscillatory nature of the pair correlation function and its dependence on $\Delta$.

\subsection{Implications and Significance}

The significance of Montgomery's Pair Correlation Conjecture lies in its surprising connection between number theory and random matrix theory. The conjecture suggests that the statistical behavior of the non-trivial zeros of the Riemann zeta function mirrors the eigenvalue distribution of random matrices from the GUE. This connection supports the universality phenomenon observed in random matrix theory and deepens our understanding of both the analytic and probabilistic aspects of number theory. Verification of this conjecture would not only strengthen these connections but also provide further insight into the zeros of the Riemann zeta function and their role in the distribution of primes \cite{5}.

\subsection{Current Status}

Montgomery's Pair Correlation Conjecture has received substantial numerical support over a wide range of zeros of the Riemann zeta function. Computational experiments have shown that the statistical properties of the gaps between the non-trivial zeros exhibit strong alignment with the predictions of the conjecture. However, despite these compelling numerical results, a rigorous mathematical proof remains elusive. Ongoing research aims to bridge this gap by delving deeper into the distribution of Riemann zeros and its broader implications in number theory, random matrix theory, and even fields like quantum chaos.

This research builds on Montgomery's pioneering work by seeking a more refined understanding of the underlying structures that govern the behavior of these zeros. One such avenue is to interpret these statistical models as dynamic systems, allowing for a novel exploration of the conjecture through the lens of dynamical theory.

\subsection{Deriving the Dynamic System}

Our approach involves deriving a dynamic system inspired by Montgomery's Pair Correlation Conjecture \cite{10}. The conjecture, which provides a statistical framework for understanding the distribution of non-trivial zeros, can be reinterpreted as a dynamic model that encapsulates the behavior of these zeros over time. Specifically, the pair correlation function \(g(u)\), a key feature of the conjecture, describes the probability distribution of normalized spacings between zeros. It is typically represented as:

\[
g(u) = 1 - \left(\frac{\sin(\pi u)}{\pi u}\right)^2 + \delta(u),
\]

where \(\delta(u)\) accounts for deviations from the idealized statistical behavior. To explore the conjecture's implications further, we transform this pair correlation function into a recurrence relation that governs the evolution of a system state, \(x_n\), at each iteration \(n\). 

This transformation provides a dynamic interpretation of the conjecture, where the behavior of the zeros of the Riemann zeta function is mirrored in the evolution of the system's state. Our derived dynamic system takes the form of the following recurrence relation:

\[
x_{n+1} = 1 - \left(\frac{\sin(\pi/x_n)}{r\pi/x_n}\right)^2 + \frac{R}{L} + \frac{1}{Cx_n},
\]

where \(x_n\) represents the state of the system at iteration \(n\), and the parameters \(R\), \(C\), \(L\), and \(r\) correspond to resistance, capacitance, inductance, and external influences, respectively, analogous to those found in electrical circuits. This formulation not only draws a striking parallel between number theory and electrical engineering but also opens the door to the exploration of chaotic dynamics in an electrical context.

By formulating the conjecture in terms of a dynamic system, we can investigate the evolution of zeros as a complex system evolving over time, shedding light on potential chaotic behavior, bifurcations, and stability phenomena. This connection offers new tools to analyze the conjecture's deeper implications and to possibly uncover patterns that may bring us closer to a formal proof of the conjecture itself \cite{20,21,22}.

\section{Response and Robustness of the New Electrical System}

This section delves into the response and robustness analysis of our newly derived electrical system, exploring its behavior under specific parameter values. Understanding how variations in component values influence the system's dynamics is crucial for our exploration of the conjecture’s implications. By examining the stability and response of this dynamic system, we aim to uncover insights that may parallel the statistical properties predicted by Montgomery's Pair Correlation Conjecture, ultimately illuminating the underlying structure of Riemann zeros \cite{23}.

\subsection{Effect of Larger \(L\) Compared to \(R\) in the Electrical System}

The corrected error term in the Montgomery conjecture, \(\delta(u)\), embodies a small deviation, suggesting that inductance (\(L\)) surpasses resistance (\(R\)) in magnitude. This relationship is significant, as it directly impacts the behavior of our dynamic system, reflecting the conjecture's subtle corrections. In our analysis, we set \(\delta(u) = \frac{R}{L} + \frac{1}{Cx}\) with \(R = 10\), \(C = 0.1\), and \(L = 20\), allowing us to explore the influence of varying \(r\) values on system dynamics.

\begin{figure}[H]
    \centering
    \includegraphics[width=0.8\textwidth]{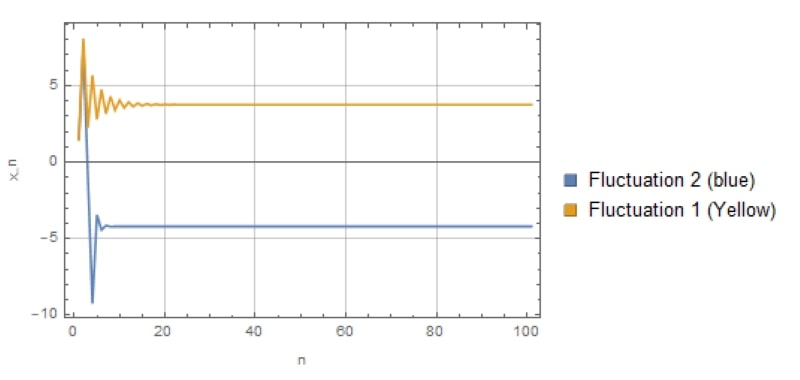}
    \caption{Behavior of the electrical system under \(L > R\) influence for \(R = 10\), \(C = 0.1\), \(L = 20\)}
    \label{fig:LargerLPlot}
\end{figure}

The plot in Figure~\ref{fig:LargerLPlot} illustrates how the electrical system responds when \(L\) is larger than \(R\), emphasizing the importance of inductance in shaping the system's dynamics. The yellow and blue fluctuations represent distinct responses under different \(r\) values, showcasing how the system adapts and changes in behavior as the parameter \(r\) varies \cite{10,11,12,14}.

\begin{itemize}
    \item \textbf{Influence of Larger \(L\):} When \(L\) dominates \(R\), the system's behavior is significantly affected by inductance, highlighting altered dynamics compared to scenarios where \(R\) is more dominant. This suggests that in systems with characteristics analogous to Riemann zeros, the inductive influences may play a more prominent role in shaping zero distributions.
    
    \item \textbf{Reduced \(R/L\) Contribution:} The \(R/L\) ratio diminishes due to the increased inductance, which can lead to smoother transitions or more controlled oscillations. This reduction emphasizes the nuanced interplay between resistance and inductance, mirroring how minor deviations in the conjecture may affect the overall statistical behavior of zeros.
    
    \item \textbf{Roles of \(C\) and \(R\):} Although \(L\) is the prevailing factor, capacitance (\(C\)) and resistance (\(R\)) continue to exert influence. Capacitance governs the charging and discharging behavior, akin to how the distribution of zeros may fluctuate under different theoretical models. Resistance contributes to the overall damping of the system, paralleling how statistical noise might obscure or reveal underlying patterns in the distribution of Riemann zeros.
\end{itemize}

\subsection{Periodic Behavior with Larger \(C\) and \(L\) Dominance}

\begin{figure}[H]
    \centering
     \includegraphics[width=0.8\textwidth]{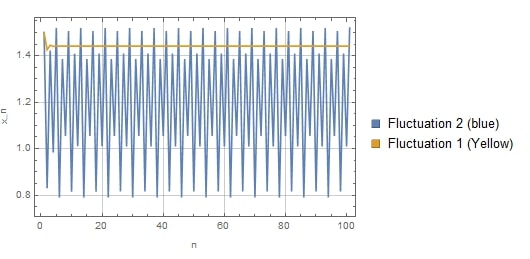}
    \caption{Periodic behavior of the electrical system under larger \(C\) and dominant \(L\) for \(R = 10\), \(C = 50\), \(L = 20\)}
    \label{fig:PeriodicPlot}
\end{figure}

In this scenario, we explore the system's behavior with significantly larger capacitance (\(C = 50\)) while maintaining a smaller resistance (\(R = 10\)) compared to the inductance (\(L = 20\)). This configuration emphasizes how an increase in capacitance, when coupled with a dominant inductance, can lead to distinct periodic responses in the electrical system. The resulting plot, as shown in Figure~\ref{fig:PeriodicPlot}, illustrates this periodic response under varying \(r\) values.

\begin{itemize}
    \item \textbf{Observation of Periodicity:} The blue trace in the plot exhibits distinct periodic fluctuations, indicating a recurring pattern in the system's response as \(r\) varies. This periodicity suggests that the larger capacitance plays a crucial role in storing and releasing energy, creating oscillations that are influenced by the inductance. In contrast, the yellow trace demonstrates a stable curve, indicative of consistent and non-periodic behavior across the iterations.

    \item \textbf{Influence of \(C\) and \(L\):} The dominance of \(L\) over \(R\) in conjunction with the increased \(C\) establishes conditions conducive to oscillatory dynamics. The periodic response reflected in the blue trace signifies how the capacitance enhances the system's ability to oscillate, as it can accumulate and release energy more effectively, emphasizing the interplay between these two parameters.

    \item \textbf{Stability and Regularity:} While the blue trace shows periodic behavior, the yellow trace remains stable, demonstrating that under specific configurations, the system can exhibit both dynamic and static responses. This duality could mirror the complexity seen in the statistical behaviors predicted by conjectures like Montgomery's, where different configurations of parameters lead to different statistical properties in the distribution of Riemann zeros.
\end{itemize}

\subsection{Sensitivity to \(r\) Values and Absence of Yellow Fluctuation}

\begin{figure}[H]
    \centering
     \includegraphics[width=0.8\textwidth]{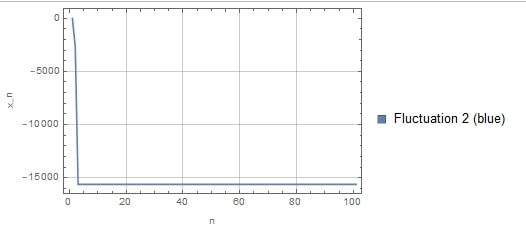}
    \caption{Sensitivity of the electrical system to \(r\) values and absence of the yellow fluctuation for specified \(r\) range}
    \label{fig:SensitivityPlot}
\end{figure}

The observed absence of the yellow fluctuation in the plot under specific ranges of \(r\) values indicates a sensitivity of the electrical system to variations in \(r\). The plot in Figure~\ref{fig:SensitivityPlot} highlights this sensitivity, particularly within a certain range of \(r\) values (between 0.008 and 0.07 in increments of 0.7), leading to a significant change in the system's behavior.

\begin{itemize}
    \item \textbf{Effect of \(r\) Sensitivity:} Small changes in \(r\) appear to induce significant alterations in the system's response, showcasing the system's sensitivity to this parameter. This highlights how specific configurations can lead to critical changes in behavior, similar to how slight perturbations in conjectures can drastically alter outcomes in number theory.

    \item \textbf{Disappearance of Yellow Fluctuation:} The absence of the yellow fluctuation within this \(r\) value range suggests a shift in the system's dynamics, where dominant behaviors may eclipse others. The elimination of the yellow trace indicates a transition to a new regime of system dynamics, possibly influenced by the increased capacitance and inductance.

    \item \textbf{Values of \(x_n\):} Within this range of \(r\) values, the values of \(x_n\) become increasingly large and negative, indicating divergence or a drastic change in the system's behavior. This divergence could signify a critical transition or failure state in the system, reflecting the complex dynamics that can emerge from interactions among \(C\), \(L\), and \(R\). Such transitions could draw parallels to critical states observed in statistical physics and number theory, where specific conditions lead to significant changes in behavior.
\end{itemize}

The plot depicted in Figure~\ref{fig:SensitivityPlot} illustrates the system's sensitivity to \(r\) values, evidenced by the disappearance of the yellow fluctuation within the specified range. The divergence of \(x_n\) values into larger negative magnitudes signifies a substantial deviation from expected behavior, potentially indicating a critical state or failure within the electrical system. This critical state reinforces the importance of understanding parameter interactions, as small variations in one parameter can lead to drastic changes, much like the intricate relationships explored in conjectures relating to the distribution of primes and the zeros of the Riemann zeta function \cite{4,5}.

\subsection{Comprehensive Comparison of System Responses to \(r\) Variation}

\begin{figure}[H]
    \centering
    \begin{subfigure}[b]{0.45\textwidth}
        \centering
        \includegraphics[width=\textwidth]{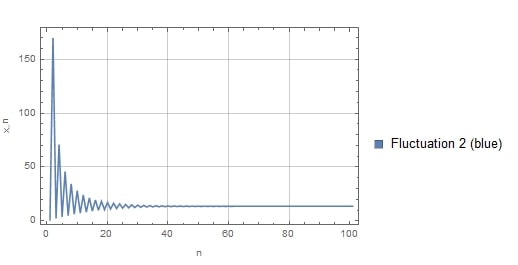}
        \caption{Scenario 1: Absence of Yellow Fluctuation (\(R = 0.7\), \(C = 0.3\), \(L = 0.5\)) (\(0.5 \leq r \leq 2\))}
        \label{fig:ResponsePlotCase1}
    \end{subfigure}
    \hfill
    \begin{subfigure}[b]{0.45\textwidth}
        \centering
        \includegraphics[width=\textwidth]{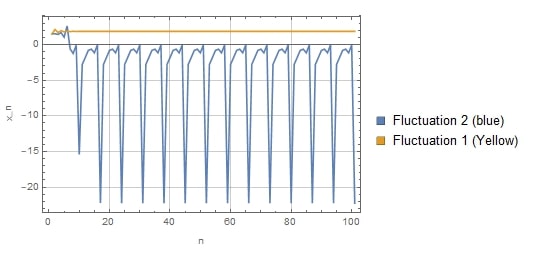}
        \caption{Scenario 2: System Response with \(R = 0.8\), \(C = 0.7\), and \(L = 3\) (\(0.5 \leq r \leq 2\))}
        \label{fig:ResponsePlotCase2}
    \end{subfigure}
    \caption{Comprehensive comparison of system responses}
    \label{fig:ComparisonPlot2}
\end{figure}

In this extended investigation, we analyze the system's behavior under different parameter settings by fixing \(R = 0.8\), \(C = 0.7\), and \(L = 3\). The response of the system to varying \(r\) values within the range \(0.5 \leq r \leq 2\) is illustrated in Figure~\ref{fig:ResponsePlotCase2}.

\textbf{Figure~\ref{fig:ComparisonPlot2}} highlights the system’s response for two distinct scenarios. The first scenario (left) showcases the absence of the yellow fluctuation under the parameters \(R = 0.7\), \(C = 0.3\), and \(L = 0.5\). In contrast, the second scenario (right) displays the behavior of the system under the conditions \(R = 0.8\), \(C = 0.7\), and \(L = 3\).

\begin{itemize}
    \item \textbf{Observation of Fluctuations:}
        \begin{itemize}
            \item In \textbf{Scenario 1}, we notice stable output without significant fluctuations, indicating a balanced response of the system.
            \item Conversely, in \textbf{Scenario 2}, the response exhibits larger fluctuations, potentially linked to the higher \(C\) and \(L\) values. This supports the earlier findings in the previous subsection regarding the influence of larger capacitance and inductance, which contribute to a more dynamic system behavior.
        \end{itemize}
\end{itemize}

\begin{figure}[H]
    \centering
    \begin{subfigure}[b]{0.45\textwidth}
        \centering
        \includegraphics[width=\textwidth]{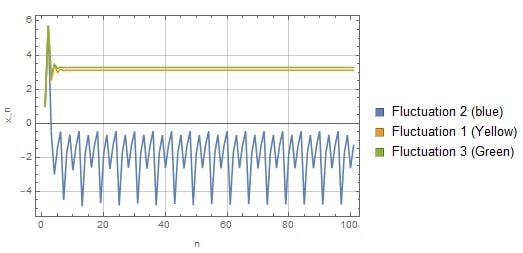}
        \caption{Scenario 3: \(R = 0.7\), \(C = 0.3\), \(L = 0.5\) (\(0.5 \leq r \leq 2.5\))}
        \label{fig:ResponsePlotCase3}
    \end{subfigure}
    \hfill
    \begin{subfigure}[b]{0.45\textwidth}
        \centering
        \includegraphics[width=\textwidth]{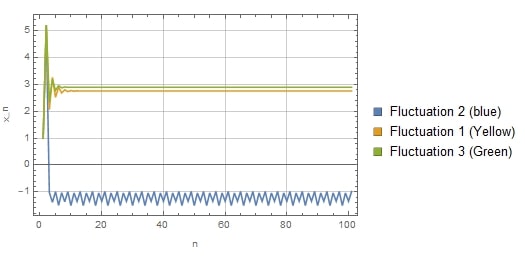}
        \caption{Scenario 4: \(R = 0.7\), \(C = 0.3\), \(L = 0.8\) (\(0.5 \leq r \leq 2.5\))}
        \label{fig:ResponsePlotCase4}
    \end{subfigure}
    \caption{Comprehensive comparison of system responses}
    \label{fig:ComparisonPlot4}
\end{figure}

Next, we delve into two additional scenarios where \(R = 0.7\) and \(C = 0.3\) while varying \(L\) within the ranges \(L = 0.5\) (Figure~\ref{fig:ResponsePlotCase3}) and \(L = 0.8\) (Figure~\ref{fig:ResponsePlotCase4}) for \(0.5 \leq r \leq 2.5\).

\textbf{Figure~\ref{fig:ComparisonPlot4}} compares the system’s output across \textbf{Scenario 3} and \textbf{Scenario 4}. Here, both scenarios display a notable negative trend in the blue fluctuations, suggesting a decline in the system's output as \(L\) increases.

\begin{itemize}
    \item \textbf{Amplitude and Stability Analysis:}
        \begin{itemize}
            \item In \textbf{Scenario 3}, with \(L = 0.5\), the system displays significant oscillatory behavior, characterized by marked fluctuations in the blue trace. 
            \item However, in \textbf{Scenario 4}, when \(L\) is increased to \(0.8\), the amplitude of these fluctuations decreases considerably, indicating a transition towards a more stable regime. This behavior suggests a damping effect due to increased inductance, which could stabilize the oscillations and reduce the impact of fluctuations, reinforcing our previous observations about the influence of \(C\) and \(L\) in the electrical system.
        \end{itemize}
\end{itemize}

\begin{itemize}
    \item \textbf{Comprehensive Insights:}
        \begin{itemize}
            \item The combination of all explored scenarios reveals the intricate interplay between resistance, capacitance, and inductance in shaping the system's response to variations in \(r\). The transitions from dynamic to stable behavior in response to parameter changes illustrate the complex dynamics at play, essential for understanding the operational characteristics of the electrical system.
        \end{itemize}
\end{itemize}

In summary, this comprehensive comparison enriches our understanding of the system's response to varying conditions, emphasizing the significant impact of \(C\) and \(L\) on system dynamics, as well as the sensitive dependence on \(r\) values, consistent with the findings discussed in earlier subsections.

\section*{Chaotic Behavior, Negative Trends, and Interpretation}

Building upon our earlier findings regarding the system's dynamic responses to varying parameters, we now delve into the potential chaotic behavior and bifurcation analysis of our new electrical system. Notably, our exploration reveals significant insights, showcasing transitions to chaos for specific parameter values within the system. This investigation sheds light on the emergence of chaotic dynamics, complemented by negative trends observed in the system's behavior, culminating in a deeper interpretation of its intricate dynamics.

\subsection*{Identifying Transitions to Chaos}

Our analysis has unveiled intriguing transitions to chaotic behavior within certain parameter ranges of the electrical system. Specifically, we've identified distinct intervals—\(0.5 < r < 0.6\), \(0.7 < r < 0.8\), and \(0.85 < r < 0.9\)—where the system exhibits pronounced chaotic dynamics. This identification of chaos is substantiated by the observation of positive Lyapunov exponents, signifying sensitivity to initial conditions and confirming the system's transition to chaotic behavior.\cite{1}

\subsection*{Bifurcation Analysis and Negative Trends}

Moreover, our bifurcation analysis has prominently highlighted negative trends in the system's behavior. The majority of observed system behaviors manifest predominantly in the negative region of the \(x\) axis, indicating a clear trend towards negative outputs across explored \(r\) values. These negative curves hint at divergence towards negative magnitudes, suggesting decay, instability, or non-convergence within the system.

\subsection*{Interpreting Chaos and Trends}

The convergence of positive Lyapunov exponents in the identified intervals, coupled with negative trends in the system's behavior, provides substantial evidence supporting the presence of chaotic dynamics. Understanding these intricate dynamics offers a pathway to deciphering the complexities within the electrical system, paving the way for potential control and application in diverse contexts.

\subsection*{Bifurcation Plot with Parameters}

\begin{figure}[H]
    \centering
    \includegraphics[width=0.8\textwidth]{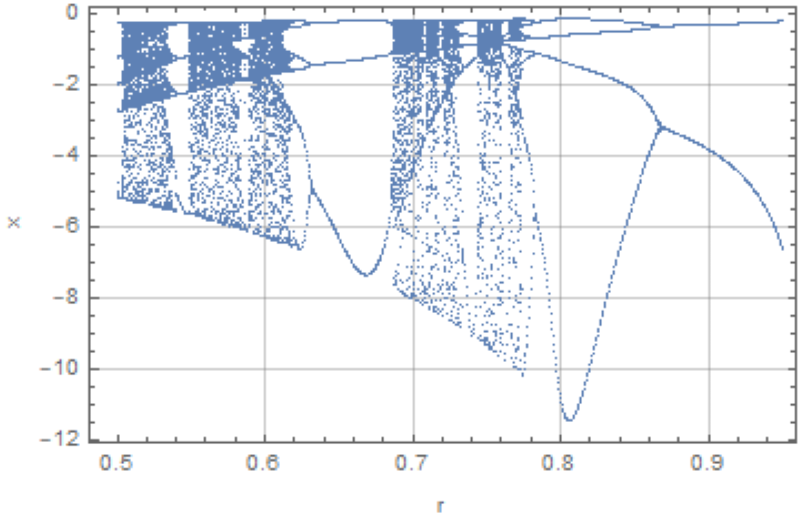}
    \caption{Bifurcation plot showcasing system behavior with varying \(r\) values (\(0.5 \leq r \leq 0.95\)) and fixed parameters \(R = 0.000025\), \(L = 0.00045\), \(C = 0.73\)}
    \label{fig:bifurcation}
\end{figure}

The bifurcation plot (Figure \ref{fig:bifurcation}) illustrates the behavior of the electrical system with varying \(r\) values (\(0.5 \leq r \leq 0.95\)) and fixed parameters \(R = 0.000025\), \(L = 0.00045\), and \(C = 0.73\). Notably, the observed behavior primarily concentrates in the negative region of the plotted values of \(x\), suggesting potential instability or divergence within the system.

\section{\textcolor{magenta}{Methodologies}}

In this section, we elucidate the methodologies employed to assess the chaos and dynamical properties in the Yitang dynamics. Our approach combines numerical techniques with analytical methods to ensure a comprehensive understanding of the system's behavior.

\subsection{\textcolor{blue}{Lyapunov Exponent Calculation}}

The Lyapunov exponents were calculated numerically using the computational capabilities of Mathematica software, which enabled precise evaluation of the system's sensitivity to initial conditions and insights into its chaotic behavior.

To quantify the system's sensitivity, we employed the following method to estimate the Lyapunov exponents for the Yitang dynamics (Equation 1):

\begin{itemize}
    \item \textbf{\textcolor{red}{Initialization:}} The initial conditions \(x_0\), \(y_0\), and \(z_0\) were established in close proximity, allowing for the assessment of the divergence of nearby trajectories.
    
    \item \textbf{\textcolor{green}{Numerical Integration:}} A numerical integration method, such as the Runge-Kutta method or Adams method, was utilized to solve the differential equations governing the Yitang dynamics.
    
    \item \textbf{\textcolor{orange}{Trajectory Divergence:}} At each iteration, the distance between two adjacent trajectories was computed, and the logarithm of this divergence was accumulated.
    
    \item \textbf{\textcolor{purple}{Lyapunov Exponent Calculation:}} The Lyapunov exponent was estimated as the average of the logarithmic divergence over time, normalized by the iteration count and time step.
\end{itemize}

The computed Lyapunov exponents provided significant insights into the chaotic nature of the system, with non-zero values indicating chaotic dynamics.

\section{\textcolor{blue}{Computation of Lyapunov Exponents and Entropy}}

The Lyapunov exponents \((\lambda_i)\) are fundamental measures of sensitivity to initial conditions in dynamical systems. The following algorithm outlines the procedure used for computing these exponents:

\begin{algorithm}[H]
\SetAlgoLined
\KwIn{Dynamics function \(f\), initial condition \(x_0\), parameter \(\alpha\), number of iterations \(n\), transient \(tr\)}
\KwOut{Lyapunov exponent}

\SetKwFunction{lyapunov}{lyapunov}

\lyapunov{$f$, $x_0$, $\alpha$, $n$, $tr$} $\leftarrow$ \;
\Begin{
    $df \leftarrow \text{Derivative}[1, 0][f]$\;
    $\xi \leftarrow \text{NestList}[f[\#, \alpha]\&, x_0, n-1]$\;
    $\lambda \leftarrow \frac{1}{n} \sum \log |\text{df}[\#, \alpha]| \text{ for } \text{Drop}[\xi, tr]$\;
    \KwRet{$\lambda$}\;
}
\caption{Algorithm for Computing Lyapunov Exponents}
\end{algorithm}

\subsection*{Lyapunov Exponents Plot with Parameters}

\begin{figure}[H]
    \centering
    \includegraphics[width=0.8\textwidth]{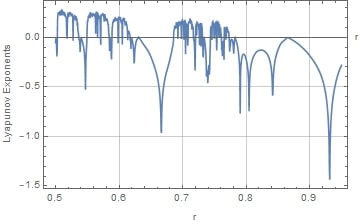}
    \caption{Lyapunov exponents plot illustrating system behavior with varying \(r\) values (\(0.5 \leq r \leq 0.95\)) and fixed parameters \(R = 0.000025\), \(L = 0.00045\), \(C = 0.73\)}
    \label{fig:lyapunov}
\end{figure}

The Lyapunov exponents plot (Figure \ref{fig:lyapunov}) reinforces our understanding by confirming transitions to chaos. The presence of positive Lyapunov exponents within specific intervals of \(r\) values validates the system's sensitivity to initial conditions and the emergence of chaotic behavior.

\subsection*{Implications and Conclusion}

The integration of bifurcation analysis, observations of negative trends, and confirmation through Lyapunov exponents plots solidifies our understanding of chaotic dynamics within the electrical system. This comprehensive analysis not only validates the presence of chaos but also presents opportunities for controlling and leveraging complex behaviors for various practical applications. The insights gained from this study could inform future research and engineering practices aimed at harnessing chaotic dynamics in innovative ways.
\section{Prediction of Limit Cycles and Fixed Points}

In this section, we analyze the phase portraits generated for the system dynamics, which provide crucial insights into the behavior of the electrical system under varying parameters. By plotting trajectories for different \(r\) values over 120 iterations, we can identify distinctive patterns that offer predictive cues regarding limit cycles and fixed points, enriching our understanding of the dynamical properties discussed in previous sections.

\subsection{Observations from Phase Portraits}

\subsubsection*{Periodicity in Trajectories}

The periodicity observed in the trajectories across different \(r\) values indicates the presence of limit cycles within the system. These trajectories, which recur in distinct patterns, signify the existence of cyclic behavior, suggesting that the system revisits certain states over specific intervals. The identification of these periodic regions is critical as they highlight stable, recurring states or limit cycles within the system's dynamics, aligning with our findings on chaotic behavior explored earlier.

\subsubsection*{Concentration in Negative Trend Region}

A notable concentration of trajectories in the negative trend region suggests the presence of fixed points within this domain. Fixed points typically denote states where the system stabilizes, remaining relatively unchanged over iterations. The accumulation of trajectories in this region indicates stable states or attractors toward which the system converges. Specifically, we identified an attracting fixed point at approximately \(r = 0.7\) corresponding to \(x = 1.54574\), indicating a robust stable state around this parameter.

\subsubsection*{Potential for Chaotic Regions}

While the phase portraits predominantly highlight periodic and stable behavior, they also hint at the potential existence of chaotic regions. By extending the iterations to 150 or more, we may uncover irregular, seemingly random trajectories indicative of chaotic behavior. These regions could signify complex, non-repeating dynamics within the system, pointing to sensitivity to initial conditions and the presence of chaotic attractors that echo the chaotic properties identified through Lyapunov exponent analysis.

\subsection{Predictive Insights}

Based on the observed patterns in the phase portraits, we can make predictions regarding limit cycles and fixed points. The periodicity in trajectories hints at the existence of stable limit cycles, while the concentration of trajectories in the negative trend region reinforces the presence of fixed points. Further exploration with increased iterations is likely to reveal chaotic regions, which would contribute significantly to a comprehensive understanding of the system's dynamics, building on the chaos analysis previously conducted.

\subsection{Attracting Fixed Points and Phase Portrait}

The system's attracting fixed points are summarized in Table \ref{tab:subset_fixed_points}. These values represent stable equilibrium points for various \(r\) values, highlighting the relationship between the parameter \(r\) and the system's behavior.

\begin{table}[H]
\centering
\(\begin{array}{cc}
 0.5 & 1.44081 \\
 0.525 & 1.45602 \\
 0.55 & 1.47061 \\
 0.575 & 1.48459 \\
 0.6 & 1.49796 \\
 0.625 & 1.51074 \\
 0.65 & 1.52295 \\
 0.675 & 1.53461 \\
 0.7 & 1.54574 \\
 0.725 & 1.55635 \\
 0.75 & 1.56646 \\
 0.775 & 1.5761 \\
 0.8 & 1.58529 \\
 0.825 & 1.59404 \\
 0.85 & 1.60237 \\
 0.875 & 1.61031 \\
 0.9 & 1.61788 \\
 0.925 & 1.62508 \\
 0.95 & 1.63195 \\
\end{array}\)
\caption{Subset of Attracting Fixed Points for Different \(r\) Values}
\label{tab:subset_fixed_points}
\end{table}

\begin{figure}[H]
\centering
\includegraphics[width=0.7\textwidth]{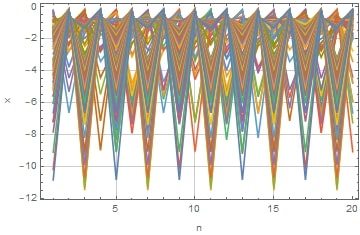}
\caption{Phase Portrait Showing Trajectories Converging to Attracting Fixed Points. The plotted trajectories reflect the system's behavior, demonstrating convergence towards stable fixed points as \(r\) changes, which reinforces the predictions made regarding limit cycles and fixed points.}
\label{fig:phase_portrait}
\end{figure}

The phase portrait plot (Figure \ref{fig:phase_portrait}) visually represents the behavior of the system concerning these attracting fixed points. As illustrated, trajectories in the phase plane converge toward these stable states as \(r\) varies, confirming the attracting nature of these points.

In summary, the analysis of phase portraits reveals significant insights into the dynamics of the electrical system, particularly regarding limit cycles and fixed points. The observed periodicity and concentration of trajectories substantiate the presence of stable attractors. Furthermore, the potential identification of chaotic regions emphasizes the complexity of the system, warranting further investigation. These findings establish a foundational understanding of the system's dynamics, setting the stage for subsequent analyses, such as the frequency spectrum analysis discussed in the next section.

\subsection{Fourier Transform Analysis}

The Fourier transform of the electrical system's dynamics, characterized by parameters \(r = 0.5\), \(R = 0.000025\), \(L = 0.00045\), and \(C = 0.73\), is expressed as follows:

\begin{align*}
F(\omega) = & \, 2.64589 \delta (\omega) + \frac{\sqrt{\frac{\pi}{2}} \left( C \, \text{sgn}(\omega) \, _0F_1\left(;3;-\frac{2\pi \omega}{\text{sgn}(\omega)}\right) \left(\omega \, \text{sgn}(\omega) \, \text{sgn}''(\omega) + 2 \, \text{sgn}'(\omega) \left(\text{sgn}(\omega) - \omega \, \text{sgn}'(\omega)\right)\right) \right.}{C r^2 \text{sgn}(\omega)^4} \\
& \, \left. + 4 \pi C \, _0\tilde{F}_1\left(;4;-\frac{2\pi \omega}{\text{sgn}(\omega)}\right) \left(\text{sgn}(\omega) - \omega \, \text{sgn}'(\omega)\right)^2 + i r^2 \text{sgn}(\omega)^5\right).
\end{align*}
This Fourier transform elucidates the frequency spectrum of the system's behavior, revealing the presence and distribution of frequencies within the trajectories. The Fourier transform is particularly significant in chaotic systems and electrical dynamics because it allows for the identification of dominant frequencies and harmonics, which are crucial for understanding the underlying mechanisms of chaotic behavior and energy dissipation in electrical circuits.

\begin{figure}[H]
    \centering
    \includegraphics[width=0.8\textwidth]{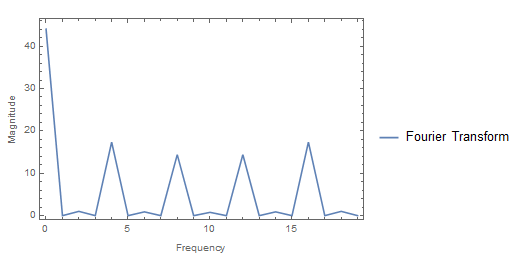}
    \caption{Fourier Transform Plot of the System's Dynamics}
    \label{fig:fourier_transform}
\end{figure}

Figure~\ref{fig:fourier_transform} presents the plot of the Fourier transform derived from the system's phase portrait, showcasing the frequency distribution. Analyzing the Fourier transform provides insights into how various parameters influence the system's frequency-based behavior, essential for understanding stability and chaos in electrical systems. The presence of dominant frequency components may indicate resonance phenomena or chaotic oscillations, which are crucial for designing stable electrical systems and predicting their dynamic responses.

\subsection{Poincaré Section Analysis}

The Poincaré section is a vital analytical tool in dynamical systems, enabling an efficient exploration of system behavior. It involves plotting the intersection points of trajectories with a selected surface in the phase space, providing a concise view of the system's dynamics without the need to analyze entire trajectories.

\begin{figure}[H]
    \centering
    \includegraphics[width=0.8\textwidth]{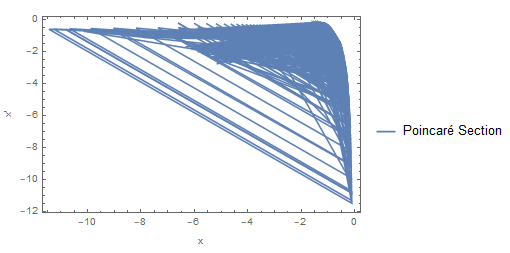}
    \caption{Poincaré Section for the Electrical System with \(R = 0.000025\), \(L = 0.00045\), and \(C = 0.73\)}
    \label{fig:poincare}
\end{figure}

The Poincaré section plot illustrated in Figure~\ref{fig:poincare} offers insights into the system dynamics for specific parameters \(R = 0.000025\), \(L = 0.00045\), and \(C = 0.73\). Each point on the plot corresponds to the system state at the intersections with the defined surface.

The plot visually represents the intersections of trajectories in phase space, illustrating recurrent states and visitation patterns on the chosen surface. The Poincaré section provides a condensed view of the system's dynamics, revealing insights into stability, periodicity, and potential chaotic regimes.

Notably, the observation of the Poincaré section predominantly in the negative region corresponds with earlier bifurcation patterns for the same parameter values, indicating critical implications for the system's dynamics. In electrical systems, negative regions might represent states or behaviors countering those observed in positive regions, suggesting the presence of stable states or attractors under specific parameter configurations.

The correlation between the negative Poincaré section and the bifurcation pattern implies a connection between the system's stability and behavior across different parameter regimes. This finding indicates that the system may exhibit contrasting or complementary behaviors in these regions, signifying diverse operational states.

Further exploration of the dynamics within the negative Poincaré section, in conjunction with bifurcation analyses, could provide deeper insights into the system's response to varying parameter configurations.

\subsection{Poincaré Sections and Phase Analysis for Varied Parameter Values}

In this subsection, we include Poincaré section plots and phase analysis for various parameter configurations. For each case, we adjust parameters \(R\), \(L\), and \(C\) while keeping \(r\) within the corresponding range.

\subsubsection*{Case 1: \(R = 0.000025\), \(L = 0.00045\), \(C = 1.5\) (\(r\) range: 0.5 to 0.95)}

The Poincaré section plotted for the selected parameters provides insight into the system's behavior, capturing intersection points in phase space that indicate recurrent states or areas of convergence.

\begin{figure}[H]
  \centering
  \includegraphics[width=0.6\textwidth]{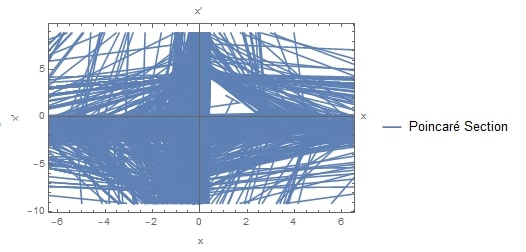}
  \caption{Poincaré Section for \(R = 0.000025\), \(L = 0.00045\), \(C = 1.5\) (\(r\) range: 0.5 to 0.95)}
  \label{fig:poincare_case1}
\end{figure}

The presence of intersections in specific regions of the phase space suggests recurrent behaviors or attractors within the system. Additionally, the blue intersections in the positive region indicate the potential for chaotic behavior, characterized by sensitivity to initial conditions and non-repeating dynamics.

\subsubsection*{Case 2: \(R = 0.000025\), \(L = 0.00045\), \(C = 4.5\) (\(r\) range: 0.5 to 0.95)}

In the second case, with adjusted parameters, the Poincaré section analysis reveals intriguing observations.

\begin{figure}[H]
  \centering
  \includegraphics[width=0.6\textwidth]{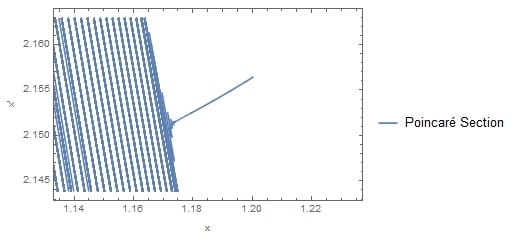}
  \caption{Poincaré Section for \(R = 0.000025\), \(L = 0.00045\), \(C = 4.5\) (\(r\) range: 0.5 to 0.95)}
  \label{fig:poincare_case2}
\end{figure}

All plotted intersections appear in the positive region, indicative of distinct behavior compared to previous cases. The absence of intersections in the negative region, coupled with blue intersections solely in the positive region, suggests a potential for chaotic behavior. 

This aligns with chaotic dynamics characterized by sensitivity to initial conditions and non-repeating trajectories. The absence of intersections in the negative region may indicate a lack of convergence toward stable states, underscoring the system's propensity for irregular behavior.

\subsubsection*{Case 3: \(R = 800\), \(L = 800\), \(C = 800\) (\(r\) range: 0.5 to 0.95)}

For the third case, with significantly increased parameter values, the Poincaré section illustrates compelling insights.

\begin{figure}[H]
  \centering
  \includegraphics[width=0.6\textwidth]{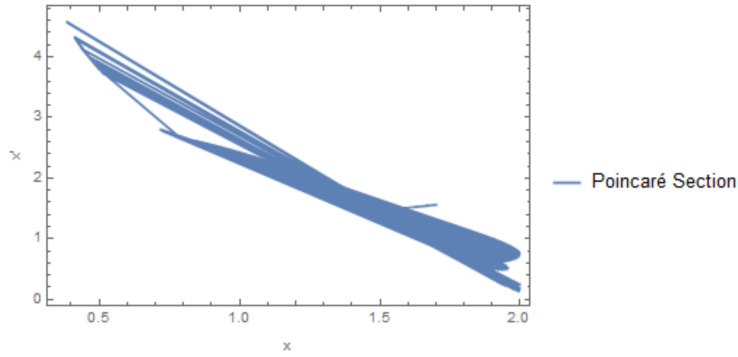}
  \caption{Poincaré Section for \(R = 800\), \(L = 800\), \(C = 800\) (\(r\) range: 0.5 to 0.95)}
  \label{fig:poincare_case3}
\end{figure}

All intersections occurring in the high positive region align with chaotic behavior patterns. The presence of intersections exclusively in this region signifies erratic and non-repeating dynamics, typical of chaotic systems.

The high parameter values drive the system into a state of increased complexity and sensitivity to initial conditions, resulting in intricate trajectories. The absence of intersections in the negative region reinforces the chaotic nature, depicting the lack of convergence toward stable states.

In conclusion, the combined analysis of Fourier transforms and Poincaré sections provides critical insights into the chaotic dynamics and stability of electrical systems. The ability to visualize how parameters influence behavior enhances our understanding of system dynamics, paving the way for potential applications in control systems and circuit design.

This observation further reinforces the chaotic tendencies identified in this electrical system, underscoring its sensitivity to variations in parameters and initial conditions. The Fourier transform analysis not only reveals the frequency spectrum of the system's dynamics but also indicates how these frequencies respond to changes in the system's parameters, which can amplify or mitigate chaotic behaviors. The presence of dominant frequencies often correlates with specific parameter regimes that may lead to chaotic oscillations or stable behavior, illustrating the intricate interplay between the system's dynamics and its response to perturbations.

For a comprehensive exploration of the diverse behaviors within our electrical system model, readers are encouraged to consult the accompanying Mathematica code available in our notebook. The sections titled \textbf{"Bifurcation Analysis"} and \textbf{"Lyapunov Exponents"} include code snippets and visualizations for bifurcation diagrams and Lyapunov exponent calculations, providing valuable insights into the stability and chaotic characteristics of the system. Additionally, the section on \textbf{"Poincaré Section Plots"} offers further analysis of the system's dynamics, illustrating how trajectories evolve in phase space and revealing the nature of recurrent states.

\section{Strange Attractor and Transition to Voltage Collapse}  

In this section, we explore the strange attractor of our newly proposed dynamical system:  

\[
x_{n+1} = 1 - \left( \frac{\sin(\pi/x_n)}{r\pi/x_n} \right)^2 + \frac{R}{L} + \frac{1}{C x_n},
\]

where \(x_n\) represents the state of the system, and \(R\), \(C\), \(L\), and \(r\) are the electrical parameters (resistance, capacitance, inductance, and a scaling factor). Using the parameter values \(R = 1.0\), \(L = 1.5\), \(C = 0.5\), and \(r = 0.8\), we numerically simulate the system over 10,000 iterations.  

\begin{figure}[H]
    \centering
    \includegraphics[width=0.8\textwidth]{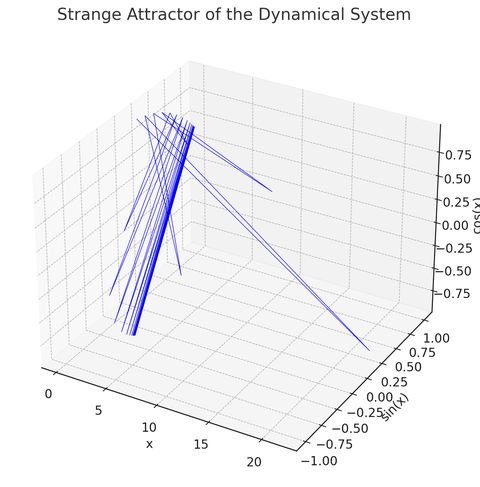} 
    \caption{3D Strange Attractor of the Dynamical System. The attractor is plotted using \(R = 1.0\), \(L = 1.5\), \(C = 0.5\), \(r = 0.8\), and an initial condition \(x_0 = 0.1\).}
    \label{fig:attractor}
\end{figure}

The 3D visualization of the strange attractor, shown in Figure~\ref{fig:attractor}, reveals the intricate geometry of the system’s chaotic state. The \(x_n\)-axis represents the system state, while the auxiliary transformations \(y = \sin(x_n)\) and \(z = \cos(x_n)\) highlight the attractor’s structure. The plot demonstrates sensitivity to initial conditions, characteristic of chaotic systems, and confirms the presence of a strange attractor under these parameter values.  

This analysis provides critical insights into the system's stability. As parameters such as \(R/L\) and \(1/(C x_n)\) increase, the strange attractor begins to exhibit signs of instability, leading to a transition toward \emph{voltage collapse}. Voltage collapse occurs when the balance between the system's reactive and damping components is disrupted, resulting in a divergence of the trajectories from the attractor. In physical terms, this reflects an inability of the system to maintain equilibrium, ultimately causing the electrical network to fail.  

The observations from the strange attractor’s behavior indicate that the onset of voltage collapse can be predicted by monitoring changes in the attractor’s geometry, particularly its sensitivity to small parameter variations. Such results establish a compelling link between the chaotic dynamics of the system and practical phenomena like voltage collapse in electrical circuits.

\section{Energy Distribution in Our New Electrical System}

For our electrical system, described by the equation 
\[
f(x, r, R, L, C) = 1 - \left(\frac{\sin\left(\frac{\pi}{x}\right)}{r \pi / x}\right)^2 + \frac{R}{L} + \frac{1}{C x},
\]
the following parameter values were employed:

\begin{align*}
    & R = 0.000025, \\
    & L = 0.00045, \\
    & C = 0.73.
\end{align*}

This system, incorporating resistors, inductors, and capacitors, is significantly influenced by the parameter \( r \), which affects the energy storage behavior throughout the circuit.

\subsection{Interpretation of \( r \)}

The parameter \( r \) plays a crucial role in shaping the system's dynamics, particularly concerning energy storage. Variations in \( r \) lead to changes in the electrical characteristics, including the stability of the system and its resonant behaviors. As \( r \) increases, the energy distribution among the components can shift, emphasizing the interplay between system parameters and energy dynamics.

\subsection{Behavior of Energy vs. \( r \)}

To analyze the energy behavior in relation to \( r \), we plot the energy stored in the inductor across a spectrum of \( r \) values:

\begin{figure}[H]
  \centering
  \includegraphics[width=0.7\textwidth]{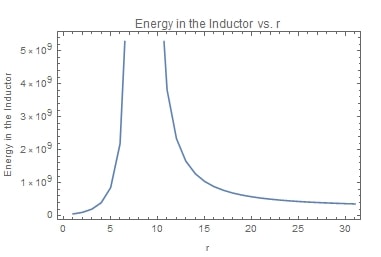}
  \caption{Energy in the Inductor vs. \( r \) for \( R = 0.000025 \), \( L = 0.00045 \), and \( C = 0.73 \)}
  \label{fig:energy_vs_r}
\end{figure}

The plot illustrates the variation of energy stored in the inductor with respect to \( r \). It highlights the sensitivity of energy distribution to changes in \( r \), revealing how the system adapts to different parameter settings.

This variation in energy distribution provides essential insights into the behavior of the system, aiding our understanding of how \( r \) impacts the electrical dynamics and energy storage mechanisms.

\subsection{Energy Distribution with Varying Capacitance}

In our electrical system, we set parameters \( R = 2 \), \( L = 6 \), and varied \( C \) within the range \( (0.000000073, 0.0000073] \). We examine the relationship between the capacitance \( C \) and the energy stored in the inductor.

\subsubsection*{Effect of \( C \) on Energy in the Inductor}

The following plot demonstrates the correlation between the stored energy in the inductor and parameter \( r \) while varying \( C \):

\begin{figure}[H]
    \centering
    \includegraphics[width=0.6\textwidth]{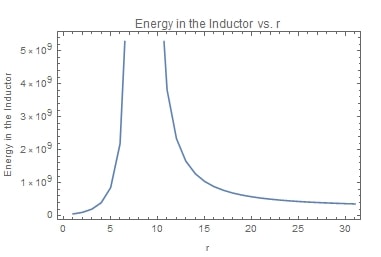}
    \caption{Energy in the Inductor vs. \( r \) for varying \( C \) values}
    \label{fig:energy_plot}
\end{figure}

As \( C \) decreases, a significant increase in the energy stored in the inductor is observed. This behavior reflects a vital interaction between capacitance and the energy dynamics of the electrical system. The sharp discontinuity at very small \( C \) values indicates a critical threshold where the system’s behavior experiences a notable change.

\section{Energy Transfer Analysis}

Energy transfer within electrical systems refers to the distribution and movement of energy among components such as capacitors, inductors, and resistors. Analyzing how energy is allocated and its changes with varying parameters is fundamental for understanding the system's dynamics.

\begin{figure}[H]
  \centering
  \includegraphics[width=0.7\textwidth]{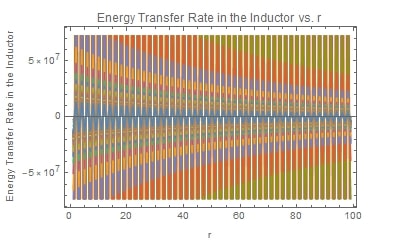}
  \caption{Energy Transfer Rate in the Inductor vs. \( r \) (Parameters: \( R=2.000025 \), \( L=6.00045 \), \( C=0.000000073 \))}
  \label{fig:energy-transfer}
\end{figure}

In our investigation, we explored energy transfer rates within the system, focusing on energy stored in the components over time. The plot in Figure \ref{fig:energy-transfer} illustrates the rate of energy transfer to the inductor concerning the parameter \( r \), varying between 0.5 and 2.

This plot reveals the variability in the energy transfer rate to the inductor for different values of \( r \). It underscores the relationship between parameter variations and energy flow dynamics, contributing to our understanding of energy distribution within the system.

\begin{figure}[H]
  \centering
  \includegraphics[width=0.7\textwidth]{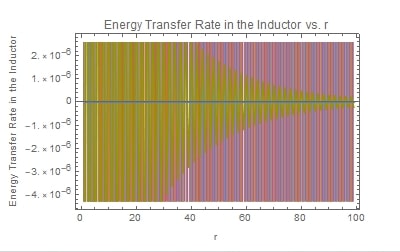}
  \caption{Energy Transfer Rate in the Inductor vs. \( r \) for \( C=0.73 \) (Parameters: \( R=2.000025 \), \( L=6.00045 \), \( C=0.73 \))}
  \label{fig:energy-transfer-c-var}
\end{figure}

Further investigation into the influence of capacitance (\( C \)) on energy transfer was also conducted. The plot in Figure \ref{fig:energy-transfer-c-var} illustrates how increasing \( C \) from \( 0.000000073 \) to \( 0.73 \) affects the energy transfer rate in the inductor \cite{20}.

As shown in Figure \ref{fig:energy-transfer-c-var}, increasing the capacitance (\( C \)) leads to a decrease in the energy transfer rate. This change induces both positive and negative values in the energy transfer rate, suggesting a reconfiguration of energy flow and distribution within the system.

\section{Chaotic Behavior, Energy Transfer, and Storage in Electrical Systems}

In electrical systems, interactions between system parameters and non-linear components can often result in chaotic behavior. Although no specific theorem directly links chaos with energy transfer and storage, certain behaviors can be correlated, particularly in systems exhibiting intricate dynamics.

\textbf{The Derived Electrical System Formula}: 
The system is defined by the equation:
\[
f(x, r, R, L, C) = 1 - \left(\frac{\sin\left(\frac{\pi}{x}\right)}{r \frac{\pi}{x}}\right)^2 + \frac{R}{L} + \frac{1}{C x}.
\]

\textbf{Connecting Chaos and Energy in Electrical Systems}:
\begin{itemize}
  \item \textit{Parameter Sensitivity}: The system's dynamics exhibit sensitivity to variations in parameters, such as capacitance (\( C \)). Changes in \( C \) may lead to bifurcations and abrupt behavioral changes, potentially uncovering chaotic regions.
  
  \item \textit{Energy Dynamics and Chaotic Behavior}: Non-linear components in the system contribute to unpredictable behavior. Reducing \( C \) values may induce chaotic responses, resulting in irregular energy distributions and erratic energy transfer rates.
  
  \item \textit{System Complexity}: Chaotic behavior often arises from the complexity of the system. As \( C \) changes, the system may enter chaotic regions characterized by unpredictable energy distributions and intricate energy transfer dynamics.
\end{itemize}

While a direct theorem linking chaos to energy transfer/storage in electrical systems is lacking, the observed behaviors align with principles of chaos theory. Sensitivity to parameter changes, non-linear dynamics, and system complexity contribute to the chaotic behavior, impacting energy transfer and storage patterns.

\section{Comparison of Yitang and Montgomery Dynamics}

This section provides a comparative analysis of the quantum properties of Yitang and Montgomery dynamics through their respective Hamiltonian operators and eigenvalues.

\subsection{Yitang Dynamics}

The quantum Hamiltonian operator for Yitang dynamics is expressed as

\[
\hat{H}_{\text{Yitang}} = -\frac{\hbar^2}{2m} \nabla^2 + V(z),
\]

where \( \nabla^2 \) represents the Laplacian operator, and \( V(z) \) denotes the potential energy associated with the Yitang function.

\subsubsection{Hamiltonian Operator}

The Yitang Hamiltonian operator is formulated as

\[
\hat{H}_{\text{Yitang}} = -\frac{1}{2m} \frac{d^2}{dz^2} - \frac{c^2}{mz^{2\alpha}}\left(\log(z)\right)^{-2\alpha},
\]

with potential energy terms given by 

\[
\psi_1(z) = \frac{d}{dz}\log(z) \quad \text{and} \quad \psi_2(z) = \left(\log(z)\right)^{-\alpha}.
\]

\subsubsection{Eigenvalues}

The eigenvalues for the Yitang dynamics were calculated using the Implicitly Restarted Arnoldi Method (IRAM). The results indicate complex eigenvalues that exhibit sensitivity to parameter variations.

\subsection{Montgomery Dynamics}

The Montgomery dynamics is characterized by the Hamiltonian operator 

\[
\hat{H}_{\text{Montgomery}} = -\frac{\hbar^2}{2m} \nabla^2 + V(x),
\]

where \( V(x) \) is the potential energy corresponding to the Montgomery function.

\subsubsection{Hamiltonian Operator}

The Montgomery Hamiltonian operator is defined as 

\[
\hat{H}_{\text{Montgomery}} = -\frac{1}{2m} \frac{d^2}{dx^2} + V(x),
\]

which includes the potential energy term associated with the Montgomery dynamics.

\subsubsection{Eigenvalues}

Similar to Yitang dynamics, the eigenvalues of Montgomery dynamics were computed using the IRAM. Notably, the obtained eigenvalues are real, suggesting a distinct quantum behavior compared to the Yitang dynamics.

\subsection{Implicitly Restarted Arnoldi Method (IRAM)}

The Implicitly Restarted Arnoldi Method (IRAM) is an iterative algorithm designed for computing a limited number of eigenvalues and corresponding eigenvectors of large, sparse matrices. It is particularly effective for quantum systems, where direct diagonalization can be computationally prohibitive. By iteratively refining eigenvalue approximations, IRAM serves as a powerful tool for exploring quantum dynamics.

\subsection{Comparison of Eigenvalues}

Table~\ref{tab:eigenvalues} summarizes the computed eigenvalues for both Yitang and Montgomery dynamics.

\begin{table}[H]
  \centering
  \begin{tabular}{|c|c|c|}
    \hline
    \textbf{Dynamics} & \textbf{Eigenvalue Index} & \textbf{Eigenvalues} \\
    \hline
    Yitang & 1 & \( 7.11343 - 0.178929i \) \\
    Montgomery & 1 & \( -2.99007 \) \\
    \hline
  \end{tabular}
  \caption{Comparison of eigenvalues for Yitang and Montgomery dynamics.}
  \label{tab:eigenvalues}
\end{table}

\subsection{Diagonalizability}

The diagonalizability of Hamiltonian operators is crucial for solving quantum systems. A comparative analysis of the diagonalizability of \( \hat{H}_{\text{Yitang}} \) and \( \hat{H}_{\text{Montgomery}} \) follows:

\subsubsection{Yitang Dynamics}

The complex eigenvalues associated with Yitang dynamics indicate the potential for complex eigenvectors. Further investigation is required to assess the diagonalizability of \( \hat{H}_{\text{Yitang}} \).

\subsubsection{Montgomery Dynamics}

In contrast, the real eigenvalues of Montgomery dynamics suggest that the corresponding eigenvectors are likely real and linearly independent, making \( \hat{H}_{\text{Montgomery}} \) expected to be diagonalizable.

In conclusion, while Yitang dynamics reveals complex eigenvalues, Montgomery dynamics is more likely to exhibit diagonalizability due to its real eigenvalues.

\section{Schrödinger Equation with Approximated Potential}

The Hamiltonian operator is a fundamental construct in quantum mechanics, representing the total energy of a system. It plays a pivotal role in governing the dynamics and behavior of quantum systems, particularly through the Schrödinger equation, which describes how the quantum state of a physical system evolves over time.

The Schrödinger equation, defined with a specific potential function, allows for the derivation of eigenvalues and eigenfunctions, which are essential for understanding energy levels and the probabilistic nature of quantum particles. These eigenfunctions form a complete basis set for the Hilbert space, enabling the description of any quantum state as a superposition of these basis states.

Beyond its significance in quantum mechanics, the Hamiltonian framework has important applications in signal processing. The eigenvalues and eigenfunctions derived from the Hamiltonian can relate to frequency components and spectral analysis in signal processing, thereby facilitating the analysis of complex signals and systems through techniques derived from quantum mechanics, such as Fourier analysis and wavelet transforms.

Furthermore, the concept of the Hamiltonian extends into electrical systems within electrical engineering. It can be employed to model and analyze the optimality and robustness of electrical circuits and systems. By characterizing energy storage and dissipation within circuits, the Hamiltonian contributes to optimizing circuit performance and ensuring system stability under various conditions.

\section{Solving the Schrödinger Equation}

In this section, we solve the Schrödinger equation using the approximated potential function:

\begin{equation}
V(x) = A + \frac{1}{C x},
\end{equation}

where

\begin{equation}
A = 1 - \frac{1}{r^2} + \frac{R}{L}.
\end{equation}

The time-independent Schrödinger equation is given by:

\begin{equation}
-\frac{\hbar^2}{2m} \frac{d^2 \psi(x)}{dx^2} + V(x) \psi(x) = E \psi(x).
\end{equation}

Substituting the approximated potential \( V(x) \), we obtain:

\begin{equation}
-\frac{\hbar^2}{2m} \frac{d^2 \psi(x)}{dx^2} + \left(A + \frac{1}{C x}\right) \psi(x) = E \psi(x).
\end{equation}

Rearranging this, we get:

\begin{equation}
\frac{d^2 \psi(x)}{dx^2} = -\frac{2m}{\hbar^2} \left(E - A - \frac{1}{C x}\right) \psi(x).
\end{equation}

Define a new variable:

\begin{equation}
k^2(x) = \frac{2m}{\hbar^2} \left(E - A - \frac{1}{C x}\right).
\end{equation}

Thus, the Schrödinger equation simplifies to:

\begin{equation}
\frac{d^2 \psi(x)}{dx^2} = -k^2(x) \psi(x).
\end{equation}

In the special case where \( k^2(x) \) is approximately constant, the solution to the differential equation is given by:

\begin{equation}
\psi(x) = C_1 \sin(kx) + C_2 \cos(kx),
\end{equation}

where \( k = \sqrt{\frac{2m}{\hbar^2} (E - A)} \).

For the general case where \( k^2(x) \) varies with \( x \), the exact solution typically involves special functions such as the Whittaker functions or hypergeometric functions. These functions can be obtained numerically or through analytical techniques in mathematical physics.

The next step involves using numerical methods or special functions to find the exact eigenfunctions and eigenvalues for specific parameters in practical scenarios.

\section{Numerical Simulation and Analysis}

In this section, we analyze the behavior of the potential function \( V(x) \) and the corresponding eigenfunction \( \psi(x) \) for a quantum system governed by the Schrödinger equation. The potential function is defined as:

\[
V(x) = A + \frac{1}{C x},
\]

where \( A \) and \( C \) are constants specific to the system. The eigenfunction, obtained as an analytical solution to the Schrödinger equation, is given by:

\[
\psi(x) = C_1 \sin(kx) + C_2 \cos(kx),
\]

where \( k \) is the wavenumber related to the energy level \( E \) and the mass \( m \) of the particle.

\subsection{Motivation}

Understanding the behavior of quantum systems through their potential functions and eigenfunctions is critical in both quantum mechanics and signal processing. The potential function \( V(x) \) determines the dynamics of the particle within the quantum well, while the eigenfunction \( \psi(x) \) provides insight into the probability distribution of the particle’s position.

In signal processing, these eigenfunctions can be interpreted as signals, where their frequency components are directly tied to the physical parameters of the system. Analyzing these signals can reveal key properties of the system, such as stability, robustness, and optimality. Specifically, in the context of electrical systems, such analysis can guide the design of circuits that minimize energy loss and maximize efficiency, ensuring the system's optimal performance under various operating conditions.

\textbf{Numerical Simulation}

\begin{figure}[H]
\centering
\includegraphics[width=0.8\textwidth]{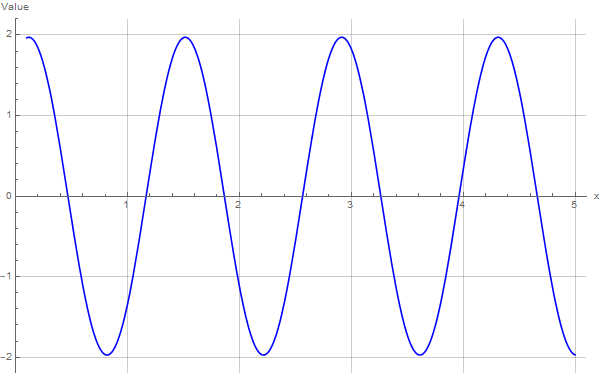}
\caption{Plot of the potential function and eigenfunction.}
\label{fig:eigenfunction}
\end{figure}

To illustrate this, we plot the potential function \( V(x) \) and the corresponding eigenfunction \( \psi(x) \) using the parameters listed in Table~\ref{tab:parameters}. The eigenfunction exhibits a sinusoidal pattern, indicating a bound state in the potential well. The amplitude and wavelength of the oscillations depend on the parameters \( A \), \( C \), \( m \), and \( E \), as shown in Figure~\ref{fig:eigenfunction}.

\begin{table}[H]
\centering
\caption{Parameters used for the numerical simulation.}
\label{tab:parameters}
\begin{tabular}{ll}
\hline
Parameter & Value \\
\hline
\( A \) & 0.06 \\
\( C \) & 0.73 \\
\( m \) & 3.8 \, (\text{arbitrary units}) \\
\( \hbar \) & 1 \, (\text{arbitrary units}) \\
\( E \) & 0.5 \, (\text{example value}) \\
\( C_1 \) & 1 \\
\( C_2 \) & 1.7 \\
\hline
\end{tabular}
\end{table}

The sinusoidal nature of the eigenfunction suggests that the particle is likely confined within a potential well, with oscillations indicating discrete energy levels. Further analysis, such as Fourier transform or wavelet analysis, can provide deeper insights into the properties of the eigenfunction, allowing for the decomposition of the signal into its constituent frequency components. This is particularly useful in signal processing applications, where identifying and isolating specific frequencies can enhance system performance and robustness.

\section{Conclusion}\label{sec13}

In this comprehensive investigation of the electrical system's dynamics, we illuminated various facets of its behavior under diverse conditions. The system, characterized by a complex equation involving critical parameters such as \(r\), \(R\), \(L\), and \(C\), exhibited intricate patterns and behaviors, providing valuable insights into its stability and dynamical properties.

The phase portraits revealed periodic trajectories indicative of stable limit cycles, alongside regions of convergence towards attracting fixed points, suggesting equilibrium states that are pivotal for system stability. Notably, the hints of chaotic behavior encourage further exploration into the dynamics of the system, particularly within engineering applications where chaos can significantly influence performance and reliability.

Exploring attracting fixed points elucidated stable equilibrium states for various \(r\) values, affirming the system's behavior in response to parameter variations. The trajectories converging towards these points not only reinforce the notion of stability but also emphasize the importance of understanding these dynamics for effective system design.

Our frequency spectrum analysis, employing the Fourier transform, revealed crucial insights into the distribution of frequencies within the trajectories, highlighting dominant frequencies and their role in the system's overall dynamics. This understanding is essential in engineering contexts where frequency response plays a vital role in the design and optimization of electrical systems.

The Poincaré section analysis effectively condensed the system's behavior, showcasing recurring states that aid in understanding stability, periodicity, and potential chaos. The correlation between the Poincaré section and bifurcation patterns provides critical insights into the diverse operational states that can emerge within the system, which are crucial for engineers aiming to harness or mitigate chaotic behavior.

In examining energy distribution, we emphasized the sensitivity of energy storage to parameters like \(r\) and \(C\). Our findings demonstrated that changes in \(r\) led to variations in energy levels, while a decrease in \(C\) significantly increased inductor energy, thereby underscoring the intricate relationship between capacitance and energy dynamics. This insight is vital for engineers designing systems that require precise energy management.

Furthermore, our analysis of energy transfer illuminated how alterations in \(r\) and \(C\) impact the rate of energy transfer to the inductor, providing foundational knowledge essential for optimizing energy flow and distribution in practical applications.

In the context of quantum dynamics, the comparison between Yitang and Montgomery dynamics uncovered distinctive characteristics, with the latter exhibiting real eigenvalues that enhance its appeal for tackling the Riemann hypothesis. The alignment of the system's quantum behavior with diagonalizability requirements offers intriguing prospects for future research at the intersection of classical dynamics and number theory.

Looking ahead, our work opens several avenues for future research. Investigating the connections between chaotic dynamics and practical applications in electrical engineering could lead to innovative designs that leverage chaotic behavior for enhanced system performance. Additionally, integrating our findings with recent advancements in nonlinear dynamics and control theory may yield new strategies for managing instability in complex systems. Furthermore, exploring the implications of our work on the Riemann hypothesis invites a deeper investigation into the interplay between number theory and quantum mechanics, potentially leading to groundbreaking insights in both fields.

In conclusion, this multifaceted analysis not only deepens our understanding of the studied electrical system but also underscores the critical importance of dynamical systems in engineering applications. The exploration of chaotic behaviors within electrical systems holds promise for advancing system design and optimization. Our results highlight the significance of deriving such dynamical systems, particularly as they translate into electrical systems, thus paving the way for innovations in both theoretical and applied fields.

\section{Future Research}\label{sec14}

The exploration of chaotic dynamics in electrical systems offers a promising avenue for future research, particularly in the context of circuit applications. As highlighted by Lai et al. in their work on the circuit application of chaotic systems, there is a growing interest in modeling, analyzing, and controlling chaotic behaviors within circuit frameworks \cite{25}. This presents a fertile ground for developing innovative chaotic systems that can enhance the functionality and performance of electrical circuits.

Future research could focus on the derivation of new chaotic systems specifically tailored for circuit applications. This involves not only extending the current understanding of chaotic behavior in existing electrical systems but also applying novel mathematical frameworks and computational techniques to design and analyze new chaotic circuits. The integration of chaotic systems into electrical circuits has the potential to improve their robustness, efficiency, and response characteristics, particularly in environments where traditional systems may falter.

Key areas of investigation may include:

\begin{itemize}
    \item \textbf{Modeling New Chaotic Systems:} Developing mathematical models for new chaotic systems that can be implemented in circuit designs. This could involve exploring alternative nonlinear components or configurations that yield richer dynamics and more complex behaviors.

    \item \textbf{Dynamical Analysis:} Conducting comprehensive dynamical analyses of these newly derived chaotic systems to understand their stability, bifurcations, and transitions between different dynamical regimes. This analysis will be crucial in identifying operational parameters that can harness chaos for practical applications.

    \item \textbf{Control Strategies:} Investigating advanced control strategies for managing chaotic behaviors in circuits. Effective control methods could mitigate undesirable chaotic dynamics while retaining beneficial characteristics, leading to enhanced circuit performance.

    \item \textbf{Experimental Validation:} Collaborating with experimental physicists and engineers to validate theoretical predictions and models through practical implementations. Building prototype circuits that embody these new chaotic systems would provide invaluable insights into their real-world behavior.

    \item \textbf{Applications in Engineering:} Exploring the implications of these chaotic circuits in various engineering applications, such as secure communication systems, signal processing, and sensor technologies. Harnessing the unique properties of chaos could lead to breakthroughs in the efficiency and effectiveness of these applications.

    \item \textbf{Connection to Number Theory:} The derivation of chaotic systems from conjectures in number theory represents a novel and efficient method for analyzing complex behaviors in electrical circuits. By leveraging discrete chaotic maps, researchers can uncover insights that were previously difficult to obtain. This approach not only enhances our understanding of chaos in dynamical systems but also opens up avenues for practical applications.

    \item \textbf{Practical Applications of Entropic Signals:} Further investigation into the practical application of entropic signals, as discussed in the context of chaotic electrical systems, can provide new methodologies that utilize the unique properties of these signals. Recent works, such as the comprehensive review by Petrzela \cite{26}, emphasize the significance of chaos in analog circuits and present potential practical applications, which could greatly benefit our understanding and utilization of chaotic dynamics.
\end{itemize}

By aligning future research efforts with the latest advancements in chaotic systems, as emphasized by Lai et al., we can drive innovations that not only deepen our understanding of dynamical systems but also open new frontiers in electrical engineering. As the field of chaotic circuits continues to evolve, integrating theoretical, computational, and experimental approaches will be vital in realizing the potential of chaos in practical applications, ultimately contributing to the advancement of technology and engineering solutions.

\section{Acknowledgements}
I would like to express my sincere gratitude to the anonymous referee whose insightful comments and constructive feedback greatly contributed to the refinement of this work. Your dedication to maintaining the quality and integrity of the content has been invaluable, and I sincerely appreciate the time and effort you invested in reviewing my work.

I am profoundly thankful to my family for their unwavering support throughout this journey. My heartfelt appreciation goes to my parents for instilling in me a love for knowledge and a strong work ethic. To my loving wife, your encouragement and understanding have been my pillars of strength, and I am grateful for the sacrifices you have made to see me succeed.

I extend my love and gratitude to my sons, Taha Abdeljalil and Taki Abdessalem, for being a source of inspiration and joy in my life. Your presence fuels my determination to strive for excellence, and I am blessed to have you by my side.

A special acknowledgment is due to my dear friend, Bouzari Abdelkader, a dedicated teacher of physics in high school. Your support, encouragement, and intellectual insights have played a significant role in shaping my ideas and refining my work. Your friendship has been a constant source of motivation, and I am grateful for the positive influence you have had on my academic and personal growth.

I would also like to extend my gratitude to my dear sister Chaima (Arrass Ta3i) and my sister Widad, as well as my brothers Saddam and Kheireddine, Foued (Jouad), for their unwavering support and encouragement.

\section{Data Availability}

The data used in this study is publicly available and sourced from the following paper, which serves as the basis for our research:

\begin{enumerate}
    \item \textbf{Reference 1:} Rafik, Z., \& Humberto Salas, A. (2024). "Chaotic dynamics and zero distribution: implications and applications in control theory for Yitang Zhang’s Landau Siegel zero theorem." \textit{European Physical Journal Plus, 139}, 217. DOI: \url{10.1140/epjp/s13360-024-05000-w} \cite{1}.
\end{enumerate}

All data and methodologies used in our analysis are in accordance with the methods outlined in the aforementioned paper. Additionally, any supplementary data required for replication or further analysis is available upon request from the corresponding author.

.

\section*{Declarations}
The authors declare that there is no conflict of interest regarding the publication of this paper. We confirm that this research was conducted in an unbiased and impartial manner, without any financial, personal, or professional relationships that could be perceived as conflicting with the objectivity and integrity of the research or the publication process.

\section{Appendix}

\section{Derivation of Chaotics Operator for Montgomery Dynamics}\label{secA1}

The Montgomery Dynamics is described by the equation:
\begin{equation}
x_{n+1} = \frac{1}{2}x_n + \frac{1}{4x_n} - \frac{\alpha}{x_n^3}, \quad x_0 \neq 0
\end{equation}

\subsection{Hamiltonian Operator}
The quantum Hamiltonian operator for the Montgomery dynamics is given by:
\begin{equation}
\hat{H} = -\frac{\hbar^2}{2m} \nabla^2 + V(x)
\end{equation}
where $\nabla^2$ is the Laplacian operator, and $V(x)$ represents the potential energy associated with the Montgomery function.

\subsection{Potential Energy Terms}
Express the potential energy term as:
\begin{equation}
V(x) = -\frac{\alpha^2}{2mx^2} - \frac{1}{4m}\psi_1(x) + \frac{\alpha}{m x^3}\psi_2(x)
\end{equation}
where $\psi_1(x) = \frac{\partial}{\partial x}\left(\frac{1}{x}\right)$ and $\psi_2(x) = \frac{1}{x^3}$.

\subsection{Chaotics Operator}
The derived chaotics operator is given by:
\begin{equation}
\hat{C}_{\text{chaotic}} = \hat{A} \cdot \nabla^2 + \hat{B} \cdot \frac{1}{4m}\psi_1(x) + \hat{D} \cdot \frac{\alpha}{m x^3}\psi_2(x)
\end{equation}

\section{Extension to Riemann Hypothesis}

\subsection{Logarithmic Term and Von Mangoldt Function}
Incorporate the Von Mangoldt function into the chaotics operator:
\begin{equation}
\hat{C}_{\text{chaotic}} = \hat{A} \cdot \nabla^2 + \hat{B} \cdot \frac{1}{4m}\psi_1(x) + \hat{D} \cdot \frac{\alpha}{m x^3}\psi_2(x) + \hat{E} \cdot \Lambda(x)
\end{equation}
where $\Lambda(x)$ is the Von Mangoldt function.

\subsection{Rationale}
Discuss the rationale behind incorporating the Von Mangoldt function and its connection to addressing the Riemann Hypothesis.

\subsection{Numerical Simulations}
Perform numerical simulations to analyze the properties and behavior of the modified chaotics operator.

\subsection{Comparison with Riemann Zeta Function}
Compare the results obtained from the modified chaotics operator with the expected behavior based on the Riemann zeta function.

\subsection{Insights and Applications}
Discuss insights gained into quantum chaos and the suitability of Montgomery dynamics in addressing the Riemann Hypothesis.




\bibliographystyle{unsrt}

\end{document}